\documentclass[12pt,a4paper]{article}
\usepackage{amsmath,amsfonts,amssymb,amsthm,amscd}

\usepackage{amsmath,amsfonts,amssymb,amsthm,amscd}
\usepackage{hyperref}

\usepackage{epstopdf}

\usepackage{braket}
\usepackage{xypic}
\usepackage{amsmath}

\usepackage[american]{babel}
\usepackage{color} 
\usepackage{epsfig}
\usepackage{caption}
\usepackage{rotating}
\usepackage{setspace}
\usepackage{fancyhdr}
\usepackage{booktabs}
\usepackage{mathrsfs}
\usepackage{amsthm}
\usepackage{amsmath,amssymb}
\usepackage{enumerate}
\usepackage{tikz}
\usepackage{youngtab}
\usepackage{hyperref}
\usepackage{tensor}
\usepackage{mathabx}
\usepackage{booktabs}
\usepackage{mathrsfs}
\usepackage{subfigure}
\usepackage{amsmath,amssymb}
\usepackage{tikz-cd}

\renewcommand{\[}{\begin{equation}}
\renewcommand{\]}{\end{equation}}

\def\id{{\rm id}}

\renewcommand{\[}{\begin{equation}}
\renewcommand{\]}{\end{equation}}

\newtheorem{thm}{Theorem}[section]
\newtheorem{cor}[thm]{Corollary}
\newtheorem{lem}[thm]{Lemma}
\newtheorem{prop}[thm]{Proposition}

\theoremstyle{definition}
\newtheorem{defn}{Definition}[section]

\theoremstyle{remark}
\newtheorem{oss}{Remark}[section]

\newcommand{\ra}{\rightarrow}

\newcommand{\zz}{\mathbb{Z}}
\newcommand{\nn}{\mathbb{N}}

\newcommand{\rr}{\mathbb{R}}

\DeclareMathOperator{\tr}{tr}

\DeclareMathOperator{\Aut}{Aut}

\DeclareMathOperator{\Tr}{Tr}

\DeclareMathOperator{\supp}{supp}
\DeclareMathOperator{\SL}{SL}

\DeclareMathOperator{\diam}{diam}

\DeclareMathOperator{\Imm}{Im}

\DeclareMathOperator{\dom}{dom}

\title{Local trivializations of suspended minimal Cantor systems and the stable orbit-breaking subalgebra}

\author{Jacopo Bassi \\ \\Department of Mathematics, University of Tor Vergata, \\Via della Ricerca Scientifica 1, 00133, Roma, Italy\\ \\email: bssjcp01@uniroma2.it}

\begin{document}

\date{}
\maketitle

\begin{abstract}
\noindent It is introduced an analogue of the orbit-breaking subalgebra for the case of free flows on locally compact metric spaces, which has a natural approximate structure in terms of a fixed point and any nested sequence of central slices around this point. It is shown that in the case of minimal flows admitting a compact Cantor central slice, the resulting $C^*$-algebra is the stabilization of the Putnam orbit-breaking subalgebra associated to the induced homeomorphism on the central slice. This construction provides an alternative characterization (up to stabilization) of the orbit-breaking subalgebra introduced by Putnam for minimal homeomorphisms of Cantor spaces in terms of suspension flows associated to such dynamical systems.
\end{abstract}

\textbf{MathScinet classification:} 46L05, 46L55 (primary)

 {\textbf{Keywords}: {$C^*$-algebra, crossed product $C^*$-algebra, dynamical system, suspension flow.}}

\section{Introduction}

In his seminal paper \cite{putnam}, Putnam showed that the structure of the crossed product $C^*$-algebra associated to a minimal homeomorphism on the Cantor set is intimately related to that of a particular $AF$-subalgebra, constructed out of Rokhlyn type decompositions of a minimal Cantor system $(X,\phi)$ with respect to nested closed sets converging to a point. This is the orbit-breaking subalgebra, as defined by Putnam and it encodes the ordered $K_0$-group and the set of traces of the crossed product, through the natural embedding; it also provides a concrete tight $AF$-embedding for the ambient algebra. It is also used in order to prove $AT$-structure for the transformation group $C^*$-algebra (this means it is approximated by algebras of the form $C(S^1, M_k)$).

Such construction can be generalized to the case of minimal homeomorphisms on compact metric spaces (\cite{lp-1,lp-2}). In this case the resulting $C^*$-algebra is approximated by Recursive Sub-Homogeneous ($RSH$) algebras, i.e. iterated pullbacks of algebras of continuous matrix-valued functions on some spaces. If the space is finite-dimensional, this $C^*$-algebra has finite decomposition rank, and this was used in \cite{toms-winter} in order to obtain a classification result for the crossed product $C^*$-algebra.

The problem of determining an analogue of the orbit breaking subalgebra in the case of actions of more general groups is open and a study of such $C^*$-algebras could bring further achievements towards a better understanding of the structure of more general transformation group $C^*$-algebras (Problem 12.7.4 and Problem 12.7.5 in \cite{cp-top}).

In this work it is defined a possible analogue of the orbit-breaking subalgebra for the case of actions of the group of real numbers on a compact metric space. It is shown that in the case the space is a suspension of a Cantor set, this $C^*$-algebra is a stable $AF$-algebra and its $K_0$-group coincides with the $K_0$-group of the crossed product; in particular it is isomorphic to the stabilization of the orbit-breaking subalgebra associated to the action of the integers on the Cantor section. In this sense this construction generalizes the original construction for minimal homeomorphisms on the Cantor set. It is expected that this $C^*$-algebra contains information about the crossed product also in more general cases where the flow cannot be reconstructed from an action of the integers on a cross section, for example the horocycle flow on the quotient of $\SL(2,\rr)$ by a discrete cocompact subgroup.

\section{Local flowbox structures}
By a flow on a topological space we mean a continuous action of the additive group of real numbers $\mathbb{R}$. 

\begin{defn}
Let $X$ be a topological space and $\phi : \rr \curvearrowright X$ a flow. $\phi$ is \textit{free} if there are no periodic points. $\phi$ is \textit{minimal} if the forward orbit of every point is dense in $X$.
\end{defn}

In the following we will use part of the notation of \cite{bartels}; in particular we will make extensive use of the notion of flowboxes and central slices, that we recall
\begin{defn}[\cite{bartels} Definition 2.3, Definition 2.4]
Let $X$ be a topological space and $\phi : \rr \ra \Aut (X)$ a flow on it. A \textit{flowbox} for $(X, \phi)$ is a compact subset $B \subset X$ such that there exists a real number $l=l_B >0$, called the \textit{length} of the box $B$, with the property that for every $x \in B$ there are real numbers $a_- (x) \leq 0 \leq a_+ (x)$ and $\epsilon (x) >0$ satisfying
\begin{equation*}
\begin{split}
l		&	=	 a_+ (x) - a_- (x) ,\\
\phi_t (x)	&	\in	 B \mbox{ for } t \in [a_- (x) , a_+ (x)],\\
\phi_t (x)	&	\notin  B \mbox{ for } t \in (a_- (x) - \epsilon (x) , a_- (x)) \cup (a_+ (x) , a_+ (x) + \epsilon (x)).
\end{split}
\end{equation*}
The set
\begin{equation*}
S_B := \{ x \in B \; | \; a_- (x) + a_+ (x) = 0 \}
\end{equation*}
is called the \textit{central slice} of $B$.
\end{defn}
If $B$ is a box, then the map
\begin{equation}
\label{eqphib}
\begin{split}
\phi_B : &S_B \times [-l/2 , +l/2] \ra B\\
	& (y, t) \mapsto \phi_t (y)
\end{split}
\end{equation}
is a homeomorphism (\cite{bartels} Lemma 2.6 (iii)) and conversely, any compact $S \subset X$ admitting an $l>0$ such that the flow restricted to $S \times [-l/2,+l/2]$ is an embedding, is the central slice of a box (this was observed in \cite{hsww} Remark 7.3.). We will usually drop the subscript $B$ in (\ref{eqphib}) and write $\phi=\phi_B$ when the flowbox $B$ is clear from the context. The following Lemma is an application of the constructions contained in Lemma 2.11 and Lemma 2.16 of \cite{bartels}.

If $X=(X,d) $ is a metric space, $x \in X$ and $r>0$, we denote by $\mathcal{B}_r(x) := \{ y \in X \; | \; d(y,x) < r\}$ the ball of radius $r$ around $x$.
\begin{lem} 
\label{lem0}
Let $X$ be a locally compact metric space and $\phi : \rr \curvearrowright X$ a free flow. For every point $x \in X$ there is a sequence of boxes $\{ B_n\}$ with non-empty interior and with associated central slices $\{ S_n\}$ and lengths $\{ l_n \}$ such that:
\begin{itemize}
\item[(i)] $S_n \subset S_{n-1} \cap B_{n-1}^\circ$ for every $i \in \nn_+$,
\item[(ii)] for every $L>0$ there exists $N \in \nn$ such that $l_n > L$ for every $n>N$,
\item[(iii)] $S_n \subset \mathcal{B}_{1/n}(x)$ for every $n \in \nn_+$, 
\item[(iv)] for every $n \in \nn$ and $0<l\leq l_n$, the box $\phi(S_n \times [-l/2, +l/2])$ has non-empty interior,
\item[(v)] for every $n \in \nn$ and $0 < \eta < l_n /2$, there is a $k >n$ such that if $L_1 \leq L_2 $ are real numbers with $-l_n /2 < L_1 - \eta$, $L_2 + \eta < l_n /2$, then 
\begin{equation}
\label{eqv1}
\phi(S_k \times [L_1, L_2]) \subset (\phi(S_n \times [L_1 - \eta, L_2 + \eta]))^\circ .
\end{equation}
\end{itemize}
\end{lem}
\proof
Let $x \in X$. Let $B_0 = \phi(S_0 \times [-l/2, +l/2])$ be the box around $x$ constructed in \cite{bartels} Lemma 2.11. Note that for every $r>0$ the box $\phi((S_0 \cap \overline{\mathcal{B}_r(x)}) \times [-l/2,+l/2])$ has non-empty interior. Let $\{l_n\}_{n \in \nn}$ be a diverging sequence of positive real numbers with $l_0 =l$. It follows from the proof of Lemma 2.16 of \cite{bartels} that there is $K \in \nn$ such that for every $k >K$ the set $S_0 \cap \overline{\mathcal{B}_{1/k} (x)}$ is the central slice of a flowbox with non-empty interior of length $l_1$; since $B_0$ has non-empty interior, it is possible to choose such a $k_1>1$ in order to obtain $\overline{\mathcal{B}_{1/k_1}(x)} \subset  B_0^\circ$ and so $S_0 \cap \overline{\mathcal{B}_{1/k_1}(x)} \subset S_0 \cap B_0^\circ$. We set $S_1 := S_0 \cap \overline{\mathcal{B}_{1/k_1} (x)}$ and $B_1 := \phi(S_1 \times [-l_1/2, +l_1/2])$. Suppose now that we are given a natural number $n >1$ and a flowbox of the form $B_n= \phi(S_n \times [-l_n/2,+l_n/2])$ with $S_n =S_0 \cap \overline{\mathcal{B}_{1/k_n} (x)}$, $k_n >n$. As above, we can find $k_{n+1} > n+1$ such that $\overline{\mathcal{B}_{1/k_{n+1}} (x)} \subset B_n^\circ $ and $S_0 \cap \overline{\mathcal{B}_{1/k_{n+1}}(x)}$ is the central slice of a flowbox with non-empty interior of length $l_{n+1}$. The sequence of flowboxes so obtained satisfies (i)-(iv); we need to prove (v). For let $S_n$ be as above, $\eta >0$ be given and $L_1< L_2$ be such that $-l_n /2 < L_1 - \eta$, $L_2+\eta < l_n/2$. From the construction it follows that $S_n$ is the central slice of a flowbox of the form $\phi(S_n \times [-\eta, + \eta])$ whose interior is non-empty and contains $x$. Hence there is $k >0$ such that $\overline{\mathcal{B}_{1/k} (x)} \subset (\phi(S_n \times [-\eta, + \eta]))^\circ$. Then
\begin{equation*}
\begin{split}
\phi_{[L_1, L_2]} (\overline{\mathcal{B}_{1/k} (x)} \cap S_n) &\subset \phi_{[L_1, L_2]} (\overline{\mathcal{B}_{1/k} (x)}) \subset  \phi_{[L_1, L_2]} (\phi((S_n \times [-\eta, +\eta]))^\circ) \\&\subset (\phi_{[L_1,L_2]} (\phi(S_n \times [-\eta, +\eta])))^\circ = (\phi_{[L_1 - \eta, L_2 + \eta]} (S_n))^\circ,
\end{split}
\end{equation*}
proving (v). $\Box$\\

\begin{defn}
\label{fs}
Let $X$ be a locally compact metric space, $\phi : \rr \curvearrowright X$ a free flow and $x \in X$. A \textit{flowbox structure (relative to the flow $\phi$) around $x$} is a sequence of boxes $\mathfrak{B} = \{ B_n\}$ with non-empty interior which satisfies conditions (i)-(v) of Lemma \ref{lem0} with respect to $x$.
\end{defn}

We recall the concept of suspension flow of a minimal homeomorphism, which gives natural examples of flowbox structures.

Let $S$ be a metric space, $\Phi \in \Aut (S)$ a free homeomorphism and $\tau : S \ra \rr_+$ a strictly positive continuous function. Define an equivalence relation on $\rr \times S$ in the following way: $(t,x) \sim (s,y)$ if there is $n \in \nn_{>0}$ such that either $y=\Phi^n (x)$ and $t=\sum_{i=0}^{n-1} \tau (\Phi^i (x)) + s$ or $x=\Phi^n (y)$ and $s=\sum_{i=0}^{n-1} \tau (\Phi^i (y))+t$. For $(t,x) \in \rr \times S$ denote by $[t,x]$ the corresponding equivalence class. Let $X_{(S, \Phi, \tau)} = (\rr \times S)/\sim$. This is the \textit{suspension} of $(S,\Phi)$ associated to the function $\tau$. The space $X_{(S,\Phi , \tau)}$ is a metric space (since $S$ is, cfr. \cite{bowen} 4) and it admits a flow defined in the following way: let $[t,x] \in X_{(S,\Phi , \tau)}$ with $0\leq t < \tau (x)$ and $s \in \rr$. If $0\leq s+t< \tau(x)$, let $\phi_s ([t,x]) = [t+s,x]$; if $s+t \geq \tau (x)$, let $k \in \nn$ be uniquely given by the condition $\sum_{i=0}^{k-1} \tau(\Phi^i (x)) \leq s+t < \sum_{i=0}^k \tau(\Phi^i (x))$ and define $\phi_s([t,x]) =[t+s-\sum_{i=0}^{k-1} \tau(\Phi^i (x)), \Phi^k (x)]$; if $t+s<0$, let $k$ be uniquely given by the condition $-\sum_{i=-k}^{-1} \tau (\Phi^i (x)) \leq s+t < -\sum_{i=-k+1}^{-1} \tau (\Phi^{-1} (x))$ and define $\phi_s ([t,x])=[s+t- \sum_{i=-k}^{-1} \tau (\Phi^{-i} (x)), \Phi^{-k} (x)]$. The space $X_{(S,\Phi, 1)}$ is sometimes called the mapping torus associated to $\Phi$.

Suppose now that we are given a compact metric space $X$ endowed with a free minimal flow and that there is a central slice $S$ with non-empty interior such that the first return map $\Phi:= \Phi_S|_S : S \ra S$ is a homeomorphism. Denote by $\tau : S \ra \rr_+$ the first return time. Since for every $[t,x] \in X_{(S,\Phi, \tau)}$ there is a unique representative $(s,y) \in \rr \times X$, with $0 \leq s < \tau(y)$, there is a well defined equivariant bijection $\psi : X_{(S,\Phi, \tau)} \ra X$. Since the map $\rr \times S \ra X$ is the composition of the quotient map $p : \rr \times X \ra X_{(S, \Phi, \tau)}$ with $\psi$, it follows that $\psi$ is continuous, hence a homeomorphism. This fact will be crucial in the foregoing sections.

\begin{oss}
Let $Y$ be an infinite compact metric space and $r \in \Aut (Y)$ a free homeomorphism. Let $X$ be the mapping torus associated to $r$, considered as a compact metric space with the metric defined in \cite{bowen}, 4. By freeness, for every $x \in Y$ and diverging sequence $\{l_n\} \subset \rr_+$ there is a diverging sequence $\{k(n)\} \subset \nn$ such that the sets $S_n:= \overline{\mathcal{B }_{1/k(n)} (y)} \cap Y $ are the central slices of flowboxes of the form $B_n =\phi (S_n \times [-l_n /2, +l_n/2])$. The sequence $\mathfrak{B} =\{ B_n\}$ is a standard flowbox structure around $y$. \end{oss}

Let $X$ be a locally compact metric space endowed with a free flow $\phi : \rr \curvearrowright X$ and $\mathfrak{B}= \{B_n = \phi (S_n \times [-l_n /2, + l_n /2])\}$ a standard flowbox structure around a point $\bar{x} \in X$. For every $n \in \nn$ define
\begin{equation*}
\mathcal{A}_n (\mathfrak{B}):= \{ a \in C_c (\rr, C_c(X)) \; | \; a (s)|_{\phi_{[0,s]}(S_n)}=0 \quad \forall \; s \geq 0 \quad \mbox{and} \quad a(s)|_{\phi_{[s,0]} (S_n)} =0 \quad \forall s \leq 0 \}.
\end{equation*}
For every $n \in \nn$, $\mathcal{A}_n (\mathfrak{B})$ is a $*$-algebra when endowed with the product and involution inherited from the crossed product. Note that the dependence on the lengths of the boxes does not appear in this definition, which actually makes sense in more generality. The transversality condition on the central slices seems a reasonable assumption, since $\rr$ is $1$-dimensional; freeness is a natural requirement too, as we shall see in Proposition \ref{prop.ind1} and the forthcoming sections. Denote by $A_n (\mathfrak{B})$ the closure of $\mathcal{A}_n (\mathfrak{B})$ inside $C(X) \rtimes \rr$. Since $\mathcal{A}_n (\mathfrak{B}) \subset \mathcal{A}_{n+1} (\mathfrak{B})$ for every $n$, we also have $A_n (\mathfrak{B}) \subset A_{n+1}(\mathfrak{B})$ for every $n$. Thus we can consider the inductive limit $(C(X)\rtimes \mathbb{R})_{\bar{x}}^\rtimes (\mathfrak{B})= \lim_n A_n (\mathfrak{B}) = \overline{\bigcup_n A_n}$, which is realized as a $C^*$-subalgebra of $C(X) \rtimes \mathbb{R}$. The next result shows that the dependence on the particular flowbox structure can be removed in the definition of this $C^*$-algebra.

Consider the $*$-algebra
\begin{equation*}
\mathcal{A}_{\bar{x}} := \{  a \in C_c (\rr, C_c(X)) \; | \; a(s)(\phi_t (\bar{x})) =0 \; \forall 0\leq t \leq s, \; s \leq t \leq 0\}.
\end{equation*}

\begin{prop}
\label{prop.ind1}
Let $\mathfrak{B}$ be a standard flowbox structure around a point $\bar{x} \in X$. Then 
\begin{equation*}
(C_0(X)\rtimes \rr)^\rtimes_{\bar{x}}:=\overline{\mathcal{A}_{\bar{x}}} =(C_0(X)\rtimes \rr)^\rtimes_{\bar{x}} (\mathfrak{B}) .
\end{equation*}
\end{prop}
\proof We only need to prove that $\mathcal{A}_{\bar{x}} \subset (C(X)\rtimes \mathbb{R})_{\bar{x}}^\rtimes$. Let $a \in \mathcal{A}_{\bar{x}}$. Denote by $\supp_\mathbb{R} (a)$ the support of the continuous function $a : \mathbb{R} \ra C_0 (X)$. Let $X' \subset X$ compact, $L>0$ with $\supp_\rr (a) \subset [-L,+L]$, $\supp (a(t)) \subset X'$ for every $t \in \rr$. Let $\epsilon >0$. We can view $a$ as a continuous function on $[-L,+L] \times X'$.

Let $\tilde{\phi} : \rr \curvearrowright \rr \times X$ be the flow defined as $\tilde{\phi}_t ((s,y))= (s, \phi_t (y))$ and let $\tilde{Y}=(\bigcup_{t \in [0,+L]} (t, \phi_{[0,t]} (x))\sqcup (\bigcup_{t \in [-L,0]} (t,\phi_{[t,0]} (x))) \subset [-L,+L] \times X$. 

Consider the map $Y: [-L,+L] \ra \mathcal{P}_c (X):= \{ \mbox{compact subsets of } X\} $, $t \mapsto Y_t = \phi_{[0,t]} (x)$ for $t \geq 0$, $t \mapsto Y_t = \phi_{[t,0]} (x)$ for $-L\leq t \leq 0$. Since the function $[-L,+L] \times X \ra \rr$, $(r,x) \mapsto d(x, \phi_r (x))$ is continuous, it follows that for every $\gamma >0$ there is a $\delta >0$ such that for $s,t \in [-L,+L]$, with $|s-t| < \delta$,
\begin{equation*}
d_H (Y_s, Y_t) \leq \diam (\phi_{|s-t|} (x)) < \gamma.
\end{equation*}
Hence $Y$ is continuous. By Lemma \ref{lem.ind1} $\tilde{Y}$ is closed in $[-L,+L] \times X$ and since it is contained in the compact set $[-L,+L] \times \phi_{[-L,+L]} (x)$, it is compact. Since $a$ vanishes on $\tilde{Y}$, it follows from Lemma \ref{lem.ind2} that there is $\eta >0$ such that $|a(s)(x)| < \epsilon /8L$ for every $(s,x) \in \tilde{\phi}_{[-\eta,+\eta]} (\tilde{Y})$. In other words, $|a(s)(\phi_t (x))| < \epsilon /8L$ for every $-\eta \leq t \leq s+\eta$, $s\geq 0$ and $s-\eta \leq t \leq +\eta$, $s \leq 0$.

Since $a$ is continuous on $[-L,+L] \times X'$, it is uniformly continuous and so there is $\rho >0$ such that for every $x \in \tilde{Y}$, $| a(s)(y)| < \epsilon/4L$ for every $y \in B_\rho (x)$, $s$ as above. It follows from Lemma \ref{lem.ind3} that there is $m \in \nn$ with associated box $B_m$ of length greater than $2L+2\eta$ such that
\begin{equation*}
|a(s)|_{\phi(S_m \times [0,t+\eta])} < \epsilon /4L \qquad \mbox{ for every } 0 \leq t \leq s\leq +L ,
\end{equation*}
\begin{equation*}
|a(s)|_{\phi(S_m \times [t-\eta,0])} < \epsilon /4L \qquad \mbox{ for every } -L \leq s \leq t \leq 0.
\end{equation*}

In virtue of Lemma \ref{lem0} we can choose $n >m$ such that $\phi (S_n \times[0,t]) \subset \phi(S_m \times [-\eta, t+\eta])^\circ$ for every $t \in [0,+L]$ and $\phi (S_n \times [t,0]) \subset \phi (S_m \times [t-\eta, +\eta])^\circ$ for every $t \in [-L,0]$. Let $Y^{(n)} : [-L,+L] \ra \mathcal{P}_c (X)$, $t \mapsto Y^{(n)}_t = \phi (S_n \times [0,t])$ for $t \in [0,+L]$, $t \mapsto Y^{(n)}_t = \phi (S_n \times [t,0])$ for $t \in [-L,0]$. Arguing as above, using compactness of $S_n$, we have that for every $\gamma >0$ there is a $\delta >0$ such that for $s,t \in [-L,+L]$, with $|s-t| < \delta$,
\begin{equation*}
d_H (Y^{(n)}_s, Y^{(n)}_t) \leq \diam (\phi_{|s-t|} (S_n)) < \gamma.
\end{equation*}
Hence $Y^{(n)}$ is continuous and $\tilde{Y}^{(n)}:=(\bigcup_{t \in [0,+L]} \phi (S_n \times [0,t])) \sqcup (\bigcup_{t \in [-L,0]} \phi (S_n \times [t,0])) \subset [-L,+L] \times \phi (S_n \times [-L,+L]) \subset [-L,+L] \times X$ is compact. 

It is easy to see that the set $\tilde{U}^{(m)} := (\bigcup_{t \in [0,+L]}(t, \phi (S_m \times [-\eta, t+\eta])^\circ)) \sqcup (\bigcup_{t \in [-L,0]}(t, \phi (S_m \times [t-\eta, +\eta])^\circ)) \subset [-L,+L] \times X$ is open.

Since $\tilde{Y}^{(n)} \subset \tilde{U}^{(m)}$ and $\tilde{Y}^{(n)}$ is compact, $\tilde{U}^{(m)}$ is open, by normality of $[-L,+L] \times X$, there is a continuous $[0,1]$-valued function $f$ which vanishes on $\tilde{Y}^{(n)}$ and takes the value $1$ on $X \backslash \tilde{U}^{(m)}$. Let $b:= fa$, where the product is intended to be the commutative product on $C([-L,+L] \times X)$. Viewing $b$ as an element of $C_c (\rr, C_c(X))$, we have
\begin{equation*}
\begin{split}
\|b-a\|& \leq \|b-a\|_1 = \int_{-L}^{+L} \sup_{x \in X} | b(s)(x) - a(s)(x)| ds \\
&= \int_{0}^{+L} \sup_{x \in \phi (S_m \times [-\eta,s+\eta])} | b(s)(x) - a(s)(x)| ds \\
&+\int_{-L}^0 \sup_{x \in \phi(S_m \times [s-\eta,+\eta])} | b(s)(x) - a(s)(x)| ds\\
&\leq L \cdot 2(\epsilon /4L) + L \cdot 2 (\epsilon /4L) = \epsilon.
\end{split}
\end{equation*}
Since $b$ belongs to $\mathcal{A}_n (\mathfrak{B})$, the result follows. $\Box$\\

\section{Suspended minimal Cantor systems}

\subsection{Rokhlin towers}
If $X$ is a locally compact space, $\phi : \rr \curvearrowright X$ is a free and minimal flow and $S \subset X$ is the central slice of a flowbox with non-empty interior, we denote the first arrive time to $S$ by $\tau_S : X \ra \rr_+$, that is $\tau_S (x) = \min \{ t >0 \; | \; \phi_t (x) \in S\}$; the first return map $\Phi_S : X \ra S$ is given by $\Phi_S (x) = \phi_{\tau_S (x)} (x)$. This map is well-defined by minimality.

We focus now on the case of a free, minimal flow admitting a central slice that is a Cantor space.

\begin{lem}
\label{lema1}
Let $\phi : \rr \curvearrowright X$ be a minimal flow and $C \subset X$ be the central slice of a flowbox $B$ with non-empty interior. Suppose that $C$ is a Cantor space and $S, T, C'$ are clopen subsets of $C$ with $S, T \subset C' \subset C \cap B^\circ$. Suppose $y \in S$ is such that $\Phi_{C'} (y) \in T$. Then there is a clopen neighborhood $V \subset S$ of $y$ such that
\begin{itemize}
\item[-] the map $V \ra C'$, $y \mapsto \Phi_{C'} (y)$ has image contained in $T$ (hence $\Phi_{C'} |_V = \Phi_{T} |_V$) and it is a homeomorphism with its image,
\item[-] the map $V \ra \rr_+$, $y \mapsto \tau_{C'} (y)=\tau_T (y)$ is continuous,
\end{itemize}
\end{lem}
\proof 
In virtue of Lemma \ref{lemaa} there is $l>0$ such that $\phi(C' \times [-l/2,+l/2]) \subset B^\circ$. Let $y \in S$ satisfy the hypothesis. Let $U \subset \phi(T \times [-l/4,+l/4])$ be an open set containing $\Phi_{C'} (y)$ (such an open set exists again in virtue of Lemma \ref{lemaa}, since $\phi(T \times (-l/4,+l/4))$ is a non-empty open set in $B^\circ$). By freeness the point $y$ admits a clopen neighborhood $V \subset S$ such that $\phi (V \times [-l/2-\tau_{C'}(y),+l/2+\tau_{C'}(y)])$ is a box in $X$ and $\phi_{\tau_{C'} (y)} (V) \subset U$.

Let $p_T : \phi(T \times [-l/4,+l/4]) \ra T$ be the projection map. The map $V \ra T$ defined as $y' \mapsto p_T (\phi_{\tau_{C'}(y)} (y'))$ is continuous, being the composition of continuous maps.

Moreover, it is injective. In fact, suppose there are $y_1, y_2 \in V$ such that $p_T (\phi_{\tau_{C'}(y)} (y_1))=p_T (\phi_{\tau_{C'}(y)}(y_2))$. Since both $\phi_{\tau_{C'}(y)}(y_1)$ and $\phi_{\tau_{C'}(y)}(y_2)$ belong to $U$, there are $s_1, s_2 \in [-l/4,+l/4]$ and $\bar{y} \in T$ such that $\phi_{\tau_{C'}(y)} (y_1) = \phi_{s_1} (\bar{y})$ and $\phi_{\tau_{C'}(y)} (y_2) = \phi_{s_2} (\bar{y})$. In particular $\phi_{-s_1} (y_1) = \phi_{-s_2} (y_2)$ and so $y_1 = \phi_{s_1 - s_2} (y_2)$. But $s_1 - s_2 \in [-l/2, + l/2]$ and so $y_1=y_2$ since they both belong to $C'$.

Hence the map $p_T \circ \phi_{\tau_{C'}(y)}  |_V : V \ra T$ is a homeomorphism with its image (since $V$ is compact). 

It follows from Lemma \ref{lemb} that, up to restricting to a smaller clopen set containing $y$, we can suppose that $\phi(V \times [l/4, \tau_{C'}(y) - l/4])$ has empty intersection with $C'$. Let $y' \in V$. Then $\phi_{\tau_{C'} (y) -l/4} (y')$ can be uniquely written as $\phi_s (z)$ for $s \in [-l/2,0]$, $z \in C'$; moreover $-s+\tau_{C'} (y) -l/4 = \tau_{C'} (y')$. Since $\phi_{\tau_{C'} (y)} (y') = \phi_{l/4 +s} (z)$ and $l/4+s \in [-l/4,+l/4]$, it follows that this is the unique way to write $\phi_{\tau_{C'} (y)} (y') = \phi_t (z')$ with $t \in [-l/4,+l/4]$, $z' \in C'$. Since $\phi_{\tau_{C'} (y)} (y')$ belongs to $\phi(T \times [-l/4,+l/4])$, it follows that $z \in T$ and $p_T (\phi_{\tau_{C'} (y)} (y')) = z = \Phi_{C'} (y')$.

Hence $p_T (\phi_{\tau_{C'}(y)}(y'))=\Phi_{C'} (y') = \phi_{\tau_{C'}(y')} (y')$ for every $y' \in V$. 

Let $y_i \ra \bar{y}$ in $V$. By continuity of $\Phi_{C'}|_V$, $\phi_{\tau_{C'}(y_i)}(y_i) \ra \phi_{\tau_{C'}(\bar{y})} (\bar{y})$. Since the sequence $\{\phi_{\tau_{C'}(y_i)}(y_i)\}$ is contained in a flowbox, we can project continuously on $\rr$ and obtain $\tau_{C'}(y_i) \ra \tau_{C'}(\bar{y})$. $\Box$

If $\mathfrak{B}$ is a flowbox structure around a point $\bar{x} \in X$ with $S_1$ a Cantor space, then we can always replace $\mathfrak{B}$ with a flowbox structure $\mathfrak{B}'$ in which $S_n$ is a clopen subset of a central slice $S$ corresponding to a flowbox $B$ with non-empty interior, such that $S_n \subset S \cap B^\circ$ for every $n \in \nn$. In this situation we will say that $\mathfrak{B}'$ is a \textit{Cantor flowbox structure}.

As a consequence of Lemma \ref{lema1} we have the following
\begin{prop}
\label{prop00}
Let $X$ be a compact metric space endowed with a free and minimal flow. Suppose that $X$ admits a Cantor central slice. Then it admits a Cantor flowbox structure $\mathfrak{B} = \{ \phi (S_n \times [-l_n /2, +l_n /2])\}_{n \in \nn}$. For every $n \in \nn$, the first return map $\Phi_n:= \Phi_{S_n}|_{S_n} : S_n \ra S_n$ is a homeomorphism and the first return time $\tau_n := \tau_{S_n} |_{S_n} : S_n \ra \rr_+$ is continuous; in particular it admits strictly positive maximum and minimum. The original flow is a suspension flow over $(S_n , \Phi_n)$.
\end{prop}
\proof Let $n \in \nn$. Take $S=T=S_n$ in Lemma \ref{lema1}. Since $S_n$ is closed and the flow is minimal, $\Phi_n (S_n) = S_n$. Lemma \ref{lema1} implies that $\Phi_n$ is continuous. Since the flow is free, it follows from the definition of the map $\Phi_n$ that it is injective. $\Box$

From now on we fix a compact metric space $X$ endowed with a free minimal flow admitting a Cantor central slice $S$. We keep the notation of Proposition \ref{prop00}.

By freeness, given a Cantor flowbox structure, we can replace it with a Cantor flowbox structure satisfying $\Phi_n^{-1} (S_{n+1}) \cap S_{n+1} = \emptyset$ for every $n \in \nn$. This hypothesis will be assumed in the sequel.

As in \cite{putnam}, we consider Rokhlin type decompositions of $S_n$ relative to the embedding $S_{n+1} \subset S_n$. Namely, let $J_{n+1,n}\subset \nn_+$ be the finite set $\tau_{S_{n+1}} (S_{n+1}) -1 \subset \nn_+$; for $j \in J_{n+1,n}$ let the \textit{floors} of the Rokhlin decomposition be the sets $F_{n+1,j}^{(k)}:= \Phi_n^k (\tau_{S_{n+1}}^{-1} (j))$, $0\leq k \leq j$. For $j \in J_{n+1,n}$ the set $T_{n+1,n}^{(j)}:=\bigsqcup_{0\leq k \leq j}F_{n+1,j}^{(k)}$ is the \textit{Rokhlin tower} with height $j$. Then $S_n = \bigsqcup_{j \in J_{n+1,n}} T_{n+1,n}^{(j)} = \bigsqcup_{j\in J_{n+1,n}} \bigsqcup_{0 \leq k \leq j} F_{n+1,j}^{(k)}$.

\subsection{The $C^*$-algebra associated to a Cantor flowbox structure}
We keep the notation of the previous sections. Again, we fix a compact metric space $X$ endowed with a free minimal flow and admitting a Cantor central slice $S$. We fix a Cantor flowbox structure $\mathfrak{B}=\{ \phi (S_n \times [-l_n /2, +l_n /2])\}_{n \in \nn}$ around a point $\bar{x} \in S$. In this section we will provide a manageable description of the $C^*$-algebra $(C(X)\rtimes \rr)_{\bar{x}}^\rtimes$, which will be used in the forthcoming section in order to compute the ordered $K_0$-group of this $C^*$-algebra. In fact, this $C^*$-algebra is an inductive limit of trivial fields of compact operators over the central slices. The information about the return times appear in the connecting morphisms. In the sequel we denote by $\mathbb{K}$ the $C^*$-algebra of compact operators on a separable Hilbert space.
\begin{prop}

\label{propa2.2}
Let $\mathfrak{B}$ be a Cantor flowbox structure around a point $\bar{x} \in X$. For every $n \in \nn$ we have $A_n (\mathfrak{B}) \simeq C(S_n, \mathbb{K})$.
\end{prop}
\proof We begin by showing that for every $n \in \nn$ $\mathcal{A}_n (\mathfrak{B})$ admits a dense $*$-embedding in $C(S_n , \mathbb{K})$. Fix $n \in \nn$. For the sake of simplicity, throughout the proof, we will write $t_y$ in place of $\tau_{S_n}(y)$ for $y \in S_n$.

Fix $n \in \nn$. Let $f \in \mathcal{A}_n (\mathfrak{B})$ and define 
\begin{equation*}
K_{n,f}(y,s,t) := t_y f(t_y (t-s)) (\phi_{t_y t} (y)), \qquad y \in S_n, \; s,t \in [0,1].
\end{equation*}

$K_{n,f}$ is (pointwise) the kernel of a convolution operator $\tilde{K}_{n,f}$ belonging to $C(S_n, \mathbb{K}(L^2([0,1])))$. Explicitly, its action on an element $\xi \in L^2 (S_n, L^2([0,1]))$ is given by
\begin{equation*}
\tilde{K}_{n,f} \xi (y,t) = \int_0^1 K_{n,f} (y,s,t) \xi (y,s) ds.
\end{equation*}
Let 
\begin{equation*}
\pi_n (f) :=\tilde{K}_{n,f} .
\end{equation*}
If $s=0$, $f(t_y t)(\phi_{t_y t} (y))=0$ for every $t \in [0,1]$. If $s=1$, then $f(t_y (t-1)) (\phi_{t_yt}(y)) = f(t_y (t-1)) (\phi_{t_y(t-1)}(\phi_{t_y} (y))=0$ for every $t \in [0,1]$, since $\phi_{t_y} (y)$ belongs to $S_n$. It $t \in \{0,1\}$, then $f (t_y (t-s))(\phi_{t_yt} (y))=0$ for every $s \in [0,1]$ since in these cases $\phi_{t_yt}(y)$ belongs to $S_n$. Hence $K_{n,f}$ belongs to $C(S_n, C_0 ((0,1) \times (0,1)))$, which, viewed as a set of convolution operators, is dense in $C(S_n , \mathbb{K}(L^2 ([0,1]))$, since $C_0 ((0,1))$ is dense in $L^2([0,1])$.

First we check that $\pi_n$ is an algebraic $*$-homomorphism. Note that, for $f,g \in \mathcal{A}_n (\mathfrak{B})$, $y \in S_n$, $s, t \in [0,1]$, we have
\begin{equation*}
(K_{n,f}K_{n,g}) (y,s,t)= \int_0^1 K_{n,f} (y,r,t) K_{n,g} (y,s,r) dr = t_y^2 \int_0^1 f(t_y (t-r)) (\phi_{t_y t} (y)) g(t_y (r-s))(\phi_{t_y r} (y)) dr,
\end{equation*}
which equals $K_{n,fg} (y,s,t)$. Moreover, $\tilde{K}_{n,f}^*$ has convolution kernel $K'$ given by $K' (y,s,t)=\overline{K}_{n,f} (y,t,s)$, which equals $K_{n,f^*}(y,s,t)$.

Now we exhibit an inverse $\pi_n^{-1} :C(S_n, C_0 ((0,1) \times (0,1))) \ra \mathcal{A}_n (\mathfrak{B})$. Let $K_n$ be an element of $ C(S_n, C_0 ((0,1) \times (0,1)))$ representing a convolution operator. Define $f_{K_n} \in C ([-1,+1], C_0(S_n \times (0,1))$ by
\begin{equation*}
f_{K_n}(s) (y,t) := \begin{cases}	\frac{1}{t_y} K_n(y,t-s, t)	&	\mbox{ for } t-1\leq s \leq t, \; 0<t<1\\
						0			&	\mbox{ otherwise }
\end{cases}.
\end{equation*}
Let $\tau_{S_n} : S_n \ra \rr_+$ be the first return time and consider the homeomorphism $\Psi_n : S_n \times (0,1) \ra  X\backslash S_n$, $(y,t) \mapsto \phi(y, \tau_{S_n} (y) t)$. Then the formula
\begin{equation*}
\pi_n^{-1} (K_n) (s):= \begin{cases} (f_{K_n}(s/\tau_{S_n}) )\circ \Psi_n^{-1} & \mbox{ on } X \backslash S_n\\
					0	&	\mbox{ on } S_n
					\end{cases}
\end{equation*}
defines an element of $C_c (\rr, C_0(X \backslash S_n)) \subset C_c (\rr, C(X))$. $\pi_n^{-1}$ is multiplicative and involutive. For $y \in S_n$, $0 \leq t \leq t_y$ and $t-t_y \leq s \leq t$, it gives
\begin{equation}
\label{eqpropa2.1}
\pi_n^{-1} (K_n) (s)(\phi_t (y))=  f_{K_n}(s/t_y) ( \Psi^{-1} (\phi_t (y))= t_y^{-1}K_n (y,(t-s)/t_y, t/t_y).
\end{equation}
By direct computation one can check that $\pi_n^{-1}(K_n) \in \mathcal{A}_n (\mathfrak{B})$. Moreover, for $y \in S_n$, $s, t \in (0,1)$, we have 
\begin{equation*}
\begin{split}
K_{n,\pi_n^{-1}(K_n)}(y,s,t)&=t_y \pi_n^{-1}(K_n) (t_y (t-s)) (\phi_{t_y t} (y))= t_y f_{K_n}(t-s)(\Psi^{-1} (\phi_{t_y t}(y)))\\
				&=t_y f_{K_n} (t-s) (y,t)=t_y [\frac{1}{t_y} K_n (y, t-(t-s), t)] = K_n (y,s,t).
\end{split}
\end{equation*}
Analogously we have $\pi_n^{-1} \circ \pi_n = \id$, hence, $*$-algebraically, $\mathcal{A}_n (\mathfrak{B}) \simeq  C(S_n, C_0 ((0,1) \times (0,1))) \subset C(S_n, \mathbb{K}(L^2 ([0,1]))$. 

Now we see how to extend $\pi_n$ to the desired $*$-isomorphism.
Let $\mu|_{S_n}$ be an invariant probability measure for the action of $\mathbb{Z}$ on $S_n$ and let $\lambda$ be the Lebesgue measure on $\mathbb{R}$. Then $\mu' := \mu|_{S_n} \times \lambda$ is an invariant measure for $X$ (cfr. \cite{hk} pag. 382) and we can normalize it to a probability measure $\mu := (\mu' (X))^{-1} (\mu|_{S_n} \times \lambda)$ on $X$. As a consequence of minimality, $\mu$ is faithful, in the sense that non-empty open sets have strictly positive measure. This induces a faithful lower semicontinuous unbounded trace $\tau_\mu$ on $C(X)\rtimes \mathbb{R}$, given, for $f \in C_c (\mathbb{R}, C(X))$, by $f \mapsto \int f(0) d\mu$, which restricts to a faithful lower semicontinuous unbounded trace on $(C(X)\rtimes \mathbb{R})_{\bar{x}}^\rtimes$. On the other hand, $\mu|_{S_n}$ induces a faithful lower semicontinuous unbounded trace on $C(S_n, \mathbb{K})$ given by $\int \Tr d\mu|_{S_n}=(\mu'(X))^{-1} (\mu|_{S_n} \times \Tr)$.

For $K_n  \in C (S_n , C_0( (0,1) \times (0,1)))$ positive, we have
\begin{equation*}
\begin{split}
	\tau_\mu (\pi_n^{-1}(K_n)) 	&=  (\mu' (X))^{-1}\int_{S_n} d\mu|_{S_n} \int_0^{t_y} dt\; \eta_n^{-1}(K_n) (0) (\phi_t (y))	\\	
				&= (\mu' (X))^{-1}\int_{S_n} d\mu|_{S_n} \int_0^{t_y} dt\; t_y^{-1} K_n (y, t/t_y , t/t_y)\\
				&= (\mu' (X))^{-1}\int_{S_n} d\mu|_{S_n} \int_0^1 dt' K_n (y,t',t')= \int \Tr d\mu|_{S_n}(K_n)
\end{split}
\end{equation*}
and so $\pi_n$ preserves these traces.

Denote by $H_{\tau_\mu}$, $H_{\int \Tr d\mu|_{S_n}}$ the GNS-Hilbert spaces and by $\lambda_{\tau_\mu} : \dom (\tau_\mu) \ra H_{\tau_\mu}$, $\lambda_{\int \Tr d\mu|_{S_n}} \ra H_{\int \Tr d\mu|_{S_n}}$ the canonical linear maps. It follows from (\cite{dixmier} 6.3.6) that $\lambda_{\tau_\mu} (\mathcal{A}_n (\mathfrak{B}))$ is dense in $H_{\tau_\mu}$ and $\lambda_{\int \Tr d\mu|_{S_n}} (\pi_n (\mathcal{A}_n (\mathfrak{B})))$ is dense in $H_{\int \Tr d\mu|_{S_n}}$. Hence
\begin{equation*}
\begin{split}
\| \pi_n (f) \| &= \| \pi_{\int \Tr d\mu|_{S_n}} (\pi_n (f))\| = \sup_{\xi \in \lambda_{\int \Tr d\mu|_{S_n}}(C(S_n , C_0 ((0,1) \times (0,1))))_1} \| \pi_{\int \Tr d\mu|_{S_n}} (\pi_n (f)) \xi\| \\
&= \sup_{\delta \in \mathcal{A}_n, \; \tau_\mu (\delta \delta^* ) =1} \int \Tr d\mu|_{S_n} (\pi_n (f\delta)\pi_n (f\delta)^*) \\
&= \sup_{\delta \in \mathcal{A}_n, \; \tau_\mu (\delta \delta^* ) =1} \tau_\mu ((f\delta)(f\delta)^*) = \| \pi_{\tau_\mu} (f) \| = \| f\|.
\end{split}
\end{equation*}
It follows that $\pi_n$ extends to a (continuous) $*$-isomorphism $A_n (\mathfrak{B}) \simeq C(S_n, \mathbb{K})$. $\Box$

For every $n \in \nn$ consider the embedding $i_{n,n+1} : A_n (\mathfrak{B})\subset A_{n+1}(\mathfrak{B})$ given by definition. Denote by $(C(X) \rtimes \mathbb{R})_{\bar{x}} \simeq (C(X) \rtimes \mathbb{R})_{\bar{x}}^\rtimes$ the inductive limit $\lim_n (C(S_n, \mathbb{K}), \iota_{n,n+1})=\overline{\bigcup_n C(S_n, \mathbb{K})}$, where for every $n \in \mathbb{N}$, $\iota_{n,n+1}:=  \pi_{n+1}i_{n,n+1} \pi_n^{-1}$.
\begin{thm}
For every $n \in \mathbb{N}$, $j \in J_{n+1,n}$ and $0\leq k \leq j$ there exists a continuous field of isometries $U_{n,n+1}^{(j,k)}: F_{n+1,j}^{(0)} \ra \mathbb{B}(L^2([0,1]))$, where $\mathbb{B}(L^2([0,1]))$ is endowed with the strict topology, such that the embedding $\iota_{n,n+1} : C(S_n, \mathbb{K}) \rightarrow C(S_{n+1}, \mathbb{K})$ is given, when evaluated on $f \in C(S_n, \mathbb{K})$, by
\begin{equation*}
\iota_{n,n+1} (f) (y) = \sum_{k=0}^j (U_{n,n+1}^{(j,k)})^* (y) f(\Phi_n^k (y)) U_{n,n+1}^{(j,k)} (y)
\end{equation*}
for $y \in F_{n+1,j}^{(0)}$.
\end{thm}
\proof For $j \in J_{n+1,n}$, $1 \leq k \leq j$, denote by $\tau_{n+1,n}^{(k)} = \tau_{F_{n+1,j}^{(k)}}|_{F_{n+1,j}^{(0)}} : F_{n+1,j}^{(0)} \ra \rr_+$ the $k$-th arrive time relative to $\Phi_n^k |_{F_{n+1,j}^{(0)}}$ and $\tau_{n+i} = \tau_{S_{n+i}}$,  $i=0,1$. For $k=0$ we let $\tau_{n+1,n}^{(0)} =0$.

Let $n \in \nn$, $j \in J_{n+1,n}$ and $0 \leq k \leq j$. Given an element $K_n^{(j,k)} \in C(F_{n+1,j}^{(k)}, C_0 ((0,1) \times (0,1)))$, define $\iota_{n,n+1} (K_n^{(j,k)}) \in C(F_{n+1,j}^{(0)}, C_0 ((0,1) \times (0,1)))$ by

\begin{equation*}
\begin{split}
\iota_{n,n+1} (K_n^{(j,k)}) &(y,s,t) =\\
&\frac{\tau_{n+1} (y)}{\tau_n (\Phi_{n}^{k} (y))}K_n^{(j,k)}\left( \Phi_n^{k} (y), \frac{\tau_{n+1} (y) s - \tau_{n+1,n}^{(k)} (y)}{\tau_n (\Phi_n^{k} (y))}, \frac{\tau_{n+1}(y) t-\tau_{n+1,n}^{(k)} (y)} {\tau_n (\Phi_n^{k} (y))} \right) 
\end{split}
\end{equation*}
for $\frac{\tau_{n+1,n}^{(k)} (y)}{\tau_{n+1} (y)} \leq t, s \leq \frac{\tau_n (\Phi_n^{k} (y)) + \tau_{n+1,n}^{(k)} (y)}{\tau_{n+1} (y)}$, $y \in F_{n+1,j}^{(0)}$ and $\iota_{n,n+1} (K_n^{(j,k)})(y,s,t)=0$ otherwise.

Let now $K_n \in C(S_n, C_0((0,1) \times (0,1))$. This can be written uniquely as $K_n = \sum_{j \in J_{n+1,n}} \sum_{k=0}^j K_n^{(j,k)}$, with $K_n^{(j,k)} \in C(F_{n+1,j}^{(k)}, C_0 ((0,1) \times (0,1)))$ as above.
Define $\iota_{n,n+1} (K_n) \in  C(S_{n+1}, C_0 ((0,1) \times (0,1)))$ by
\begin{equation*}
\iota_{n,n+1} (K_n) = \sum_{j \in J_{n+1,n}} \sum_{k=0}^j \iota_{n,n+1} (K_n^{(j,k)}) .
\end{equation*}

 We have $\iota_{n,n+1} \circ \pi_n = \pi_{n+1} \circ i_{n,n+1} : \mathcal{A}_n (\mathfrak{B})\ra C(S_{n+1}, \mathbb{K})$.

We exhibit now a continuous extension of the embedding $\iota_{n,n+1}$ to an embedding $\iota_{n,n+1} : C(S_n , \mathbb{K}) \ra C(S_{n+1}, \mathbb{K})$. For, let $j \in J_{n+1,n}$ and $0 \leq k \leq j$ be given and define the continuous field of operators $U_{n,n+1}^{(j,k)}: F_{n+1,j}^{(0)} \ra \mathbb{B}(L^2([0,1]))$ (continuity is intended in the strict topology of $\mathbb{B}(L^2 ([0,1])$)) given by

\begin{equation*}
(U_{n,n+1}^{(j,k)}(y) \xi)(t) =\begin{cases}	\sqrt{\frac{\tau_{n+1} (y)}{\tau_n (\Phi_n^k (y))}} \xi \left(\frac{\tau_{n+1} (y) t - \tau_{n+1,n}^{(k)} (y)}{\tau_n (\Phi_n^k (y))}\right) & \mbox{ for } \frac{\tau_{n+1,n}^{(k)} (y)}{\tau_{n+1}(y)} \leq t \leq \frac{\tau_n (\Phi_n^k (y)) + \tau_{n+1,n}^{(k)} (y)}{\tau_{n+1} (y)},\\
\qquad \qquad 0 & \mbox{ otherwise }.
\end{cases}
\end{equation*}
Then $U_{n,n+1}^{(j,k)}(y)$ is an isometry for every $y \in F_{n+1,j}^{(0)}$ and its adjoint is given by
\begin{equation*}
((U_{n,n+1}^{(j,k)}(y))^* \xi)(t) = \sqrt{ \frac{\tau_n (\Phi_n^k (y))}{\tau_{n+1} (y)}} \xi \left( \frac{\tau_n (\Phi_n^k (y)) t + \tau_{n+1,n}^{(k)} (y)}{\tau_{n+1} (y)}\right);
\end{equation*}
the connecting morphism $\iota_{n,n+1} : C(F_{n+1,j}^{(k)}, \mathbb{K}) \ra C(F_{n+1,j}^{(0)}, \mathbb{K}) \subset C(S_{n+1}, \mathbb{K})$ reads, for $f \in C(F_{n+1,j}^{(k)}, \mathbb{K})$, $y \in F_{n+1,j}^{(0)}$,
\begin{equation*}
\iota_{n,n+1} (f)(y) =U_{n,n+1}^{(j,k)}(y) f(\Phi_n^k (y)) (U_{n,n+1}^{(j,k)}(y))^*.
\end{equation*}
Hence, writing $f \in C(S_n , \mathbb{K})$ as $f=\sum f_{j,k}$, with $f_{j,k} \in C(F_{n+1,j}^{(k)}, \mathbb{K})$, we obtain $\iota_{n,n+1} (f) = \sum \iota_{n,n+1} (f_{j,k})$. We note that if $y \in F_{n+1,j}^{(0)}$, the associated isometries $U_{n,n+1}^{(j,k)}(y)$, $k=0,...,j$ have pairwise orthogonal ranges. In particular, $\iota_{n,n+1}$ preserves the orthogonality of the $f_{j,k}$'s and defines a $*$-homomorphism from $C(S_n, \mathbb{K})$ to $C(S_{n+1}, \mathbb{K})$. $\Box$

\subsection{The $K_0$-group}

This section strongly relies on the results and ideas contained in \cite{putnam}. Our purpose here is to obtain an analogue of the exact sequence appearing in Theorem 4.1 of \cite{putnam} for our situation.

In the following we will denote $\iota_n = \iota_{n,n+1} : C(S_n, \mathbb{K}) \ra  C(S_{n+1}, \mathbb{K})$. If $\phi : A \ra B$ is a $*$-homomorphism between $C^*$-algebras we let $\phi_* : K_0 (A) \ra K_0 (B)$ be the corresponding group homomorphism. We fix a Cantor flowbox structure satisfying the hypothesis of the previous section.

Since $K_0 (C(Y))= C(Y , \mathbb{Z})$ for every totally disconnected compact space $Y$, we can write for $f \in C(S_n , \mathbb{Z})$ and $y \in F_{n+1,j}^{(0)}$, with $j \in J_{n+1,n}$ and $0 \leq k \leq j$,
\begin{equation*}
\iota_{n,*}(f)(y)=\sum_{k=0}^j (f(\Phi_{n}^{k}(y)) = \sum_{k=0}^j (\hat{\Phi}_{n,\mathbb{K}, *}^k (f))(y),
\end{equation*}
where we denote by $\hat{\Phi}_{n,\mathbb{K}} \in \Aut (C(S_n, \mathbb{K}))$ the natural extension to $C(S_n , \mathbb{K})$ of the action $\hat{\Phi}_n$ on $C(S_n)$. Under the $*$-isomorphism $C(S_n , \mathbb{K}) \simeq C(S_n) \otimes \mathbb{K}$ this is the same as $\hat{\Phi}_n \otimes \id$.

For every $n \in \nn$ we denote $D_n := C(S_n, \mathbb{K})$ and fix a $*$-isomorphism $C(S_n \rtimes_{\Phi_n} \mathbb{Z}) \otimes \mathbb{K} \simeq C(X) \rtimes \mathbb{R}$ (\cite{hsww} Corollary 9.1 and Corollary 6.7, \cite{rme}). By the Pimsner-Voiculescu exact sequence, the inclusion $\alpha_n : C(S_n) \subset C(S_n) \rtimes_{\Phi_n} \mathbb{Z}$ gives a surjection at the level of $K_0$-groups. Hence the same is true for the composition $\gamma_n : D_n =C(S_n, \mathbb{K}) \simeq C(S_n) \otimes \mathbb{K} \subset C(S_n \rtimes_{\Phi_n} \mathbb{Z}) \otimes \mathbb{K} \simeq C(X) \rtimes \mathbb{R}$ for every $n \in \nn$.

\begin{prop}
\label{propn1}
Let $n \in \nn$ and $f \in C(S_n , \mathbb{Z})= K_0 (D_n)$ be such that $f|_{\Phi_n^{-1} (S_{n+1})}=0$. Then $\iota_{n,*} (f) = \iota_{n,*} (\hat{\Phi}_{n,\mathbb{K},*} (f))$.\end{prop}
\proof 
Let $j \in J_{n+1,n}$ and $y \in F_{n+1,j}^{(0)} \subset S_{n+1}$; denote by $y_1=\Phi_n (y), y_2=\Phi_n^2 (y),..., y_i=\Phi_n^j (y)$ the points of $S_n$ where $y$ lands following the positive direction of the flow before returning to $S_{n+1}$, ordered accordingly to the natural order of the real numbers. Then
\begin{equation*}
\iota_{n,*} (f)(y)=\tr \left(\begin{array}{cccc}	f(y)	&	0		&			&	\\
								0	&	f(y_1)	&	0		&	\\
									&			&	\ddots	&	\\
									&			&		0	&	f(y_j)
									\end{array}\right)= \tr \left(\begin{array}{cccc}	f(y)	&	0		&			&	\\
								0	&	\ddots	&	0		&	\\
									&			&	f(y_{j-1})	&	\\
									&			&		0	&	0
									\end{array}\right).
\end{equation*}
On the other hand,
\begin{equation*}
\begin{split}
\iota_{n,*} (\hat{\Phi}_{n,\mathbb{K},*} (f))(y)&=\tr \left(\begin{array}{cccc}	f(\Phi_n^{-1}(y))	&	0		&			&	\\
								0	&	f(\Phi_n^{-1}(y_1))	&	0		&	\\
									&			&	\ddots	&	\\
									&			&		0	&	f(\Phi_n^{-1}(y_j))
									\end{array}\right)\\
									&= \tr \left(\begin{array}{cccc} 0	&	0		&			&	\\
								0	&	f(y)	&	0		&	\\
									&			&	\ddots	&	\\
									&			&		0	&	f(y_{j-1})
									\end{array}\right).
									\end{split}
\end{equation*}
The result follows. $\Box$

In virtue of Proposition \ref{propn1}, we can define, for every $n \in \nn$, a map $C(\Phi_n^{-1} (S_{n+1}), \mathbb{Z}) \ra K_0 (D_{n+1})$ in the following way. Let $f \in C(\Phi_n^{-1} (S_{n+1}, \mathbb{Z}))$, extend it to a continuous function $g \in C(S_n , \mathbb{Z})$. Then the map
\begin{equation*}
\beta_{n+1} (f) = \iota_{n,*} (g) - \iota_{n,*} (\hat{\Phi}_{n,\mathbb{K},*} (g)) \in K_0 (D_{n+1})
\end{equation*}
is well defined. Indeed, if $g_1, g_2 \in C(S_n, \mathbb{Z})$ are two extensions of $f$, then by Proposition \ref{propn1} $\iota_{n,*} (g_1) - \iota_{n,*} (g_2) = \iota_{n,*} (g_1 - g_2) = \iota_{n,*} (\hat{\Phi}_{n,\mathbb{K},*} (g_1-g_2)) = \iota_{n,*} (\hat{\Phi}_{n,\mathbb{K},*} (g_1)) - \iota_{n,*} (\hat{\Phi}_{n,\mathbb{K},*} (g_2))$. 
\begin{prop}
\label{propn2}
For every $n \in \nn$, there is a sequence
\begin{equation*}
0 \ra \mathbb{Z} \overset{\eta_{n+1}}{\ra} C(\Phi_n^{-1} (S_{n+1}), \mathbb{Z}) \overset{\beta_{n+1}}{\ra} K_0 (D_{n+1}) \overset{\gamma_{n+1,*}}{\ra} K_0 (C(X) \rtimes \mathbb{R}) \ra 0,
\end{equation*}
which is exact everywhere, except possibly in $C(\Phi^{-1}_n(S_{n+1}), \mathbb{Z})$, where we have $\Imm (\eta_{n+1}) \subset \ker (\beta_{n+1})$.
\end{prop}
\proof 
The map $\eta_{n+1} : \mathbb{Z} \ra C(\Phi_n^{-1} (S_{n+1}), \mathbb{Z})$ sends an integer to the corresponding constant-valued function on $\Phi_n^{-1} (S_{n+1})$. Hence exactness in $\mathbb{Z}$ is clear. It is also clear that $Im (\eta_{n+1}) \subset \ker (\beta_{n+1})$, since when applying $\beta_{n+1}$ we can choose constant extensions, which vanish under $\beta_{n+1}$.\\
We want to prove exactness in $K_0 (D_{n+1})$. For a generic choice of $i \in J_{n+1,n}$ and a point $x \in F_{n+1,i}^{(0)} \subset S_{n+1}$, we denote by $x_1=\Phi_n (x),..., x_i=\Phi_n^i (x)$ its forward orbit in $S_n \backslash S_{n+1}$ and by $x_{i+1}=\Phi_n^{i+1}(x)$ its first return to $S_{n+1}$. Let $g \in C(\Phi_{n}^{-1} (S_{n+1}), \mathbb{Z})$ and extend it to zero outside $\Phi_{n}^{-1} (S_{n+1})$; call this extension $f$. Let $j \in J_{n+1,n}$, $y  \in F_{n+1,j}^{(0)} \subset S_{n+1}$, $y=y_0, y_1 = \Phi_n (y),y_2=\Phi_n^2 (y),..., y_j=\Phi_n^j (y)$ its forward orbit, then $y= \Phi_{n} (x_i) = x_{i+1}$ for some $i \in J_{n+1,n}$, $x_i \in \Phi_n^{-1}(S_{n+1})$, with $x=x_0 \in F_{n+1,i}^{(0)}$. Note that $\Phi_{n+1}^{-1} (y) =x_0=x$. We have

\begin{equation*}
\begin{split}
(\hat{\Phi}_{n+1, \mathbb{K},*}( \iota_{n, *}(f)))(y)&=(\iota_{n, *}(f))(x)= \tr \left(\begin{array}{cccc}	f(x)	&	0		&			&	\\
								0	&	f(x_1)	&	0		&	\\
									&			&	\ddots	&	\\
									&			&		0	&	f(x_i)\end{array}\right) \\
									&=\tr \left(\begin{array}{cccc}	0	&	0		&			&	\\
								0	&	0	&	0		&	\\
									&			&	\ddots	&	\\
									&			&		0	&	f(x_i)\end{array}\right),\end{split}	
\end{equation*}

while, for the above choice of $j$,
\begin{equation*}
\iota_{n,*} (\hat{\Phi}_{n,\mathbb{K},*} (f) ) (y)=\tr \left(\begin{array}{cccc}	f(\Phi_n^{-1}(y))	&	0		&			&	\\
								0	&	f(\Phi_n^{-1}(y_1))	&	0		&	\\
									&			&	\ddots	&	\\
									&			&		0	&	f(\Phi_n^{-1}(y_j))\end{array}\right).
\end{equation*}
Note now that for every $1 \leq k \leq j$, $\Phi_n^{-1} (y_k) \in \Phi_n^{-1} (S_{n+1})$ if and only if $y_k \in S_{n+1}$; then, since $\Phi_n^{-1} (y) = x_i$, we have 
\begin{equation*}
\iota_{n,*} (\hat{\Phi}_{n,\mathbb{K},*} (f)) (y)=\tr \left(\begin{array}{cccc}	f(x_i)	&	0		&			&	\\
								0	&	0	&	0		&	\\
									&			&	\ddots	&	\\
									&			&		0	&	0\end{array}\right).
\end{equation*}
Then we obtain
\begin{equation*}
\hat{\Phi}_{n+1,\mathbb{K},*} (\iota_{n,*} (f)) = \iota_{n,*} (\hat{\Phi}_{n,\mathbb{K},*} (f)),
\end{equation*}
from which it follows that
\begin{equation*}
\gamma_{n+1,*} \iota_{n,*} (f) = \gamma_{n+1,*} \hat{\Phi}_{n+1,\mathbb{K},*} \iota_{n,*} (f) = \gamma_{n+1,*} \iota_{n,*} \hat{\Phi}_{n,\mathbb{K},*} (f).
\end{equation*}
proving that $Im (\beta_{n+1}) \subset \ker (\gamma_{n+1, *})$. In order to prove the other inclusion, we make use of the Pimsner-Voiculescu exact sequence and the commutativity of the diagram
\[ \begin{tikzcd}
K_0 (C(S_{n+1})) \arrow{rr}{\id_* - \hat{\Phi}_{n+1,*}} \arrow[swap]{d}{(\id \otimes e_{1,1})_*} && K_0 (C(S_{n+1})) \arrow{d}{(\id \otimes e_{1,1})_*} \arrow{rr}{\alpha_{n+1,*}} && K_0 (C(S_{n+1}) \rtimes_{\hat{\Phi}_{n+1}} \mathbb{Z}) \arrow{d}{(\id \otimes e_{1,1})_*} \\%
K_0 (C(S_{n+1}) \otimes \mathbb{K}) \arrow{rr}{\id_* - (\hat{\Phi}_{n+1} \otimes \id)_*}&& K_0 (C(S_{n+1}) \otimes \mathbb{K}) \arrow{rr}{(\alpha_{n+1} \otimes \id)_*} && K_0 ((C(S_{n+1}) \rtimes_{\hat{\Phi}_{n+1}} \mathbb{Z})\otimes \mathbb{K}),
\end{tikzcd}
\]
which give $\ker (\gamma_{n+1, *}) = \Imm (\id_* - \hat{\Phi}_{n+1, \mathbb{K},*})$, where $\id_* - \hat{\Phi}_{n+1,\mathbb{K}, *} : C(S_{n+1}, \mathbb{Z}) \ra C(S_{n+1}, \mathbb{Z})$. Hence we only need to observe that any element in $C(S_{n+1}, \mathbb{Z})$ is the image under $\iota_{n,*}$ of a $\mathbb{Z}$-valued function on $S_n$ which vanishes outside $\Phi_{n}^{-1} (S_{n+1})$, since, as just seen, for such functions we have $\hat{\Phi}_{n+1, \mathbb{K},*} (\iota_{n,*} (f)) = \iota_{n,*}(\hat{\Phi}_{n,\mathbb{K},*} (f))$. In order to see this, let $f \in C(S_{n+1} , \mathbb{Z})$ and consider the homeomorphisms $\Phi_{n}^{i} : F_{n+1,i}^{(0)} \ra \Phi_n^{-1} (S_{n+1})$, for $i \in J_{n+1,n}$. Hence $\Phi_n^{-1} (S_{n+1}) = \bigsqcup_{i \in J_{n+1,n}} \Phi_{,n}^{i} (F_{n+1,i}^{(0)})$. Define the function $g \in C(S_n , \mathbb{Z})$ given by $g|_{\Phi_n^{-1}(S_{n+1})} = \sum_{i\in J_{n+1,n}} \chi_{\Phi_{n}^{i} (F_{n+1,i}^{(0)})} (f \circ  (\Phi_{n}^{i})^{-1} )$ and $g|_{S_n \backslash \Phi_{n}^{-1} (S_{n+1})} =0$. Then $\iota_{n,*} (g) =f$. $\Box$\\

For every $n \in \nn$ let $e_{n} : C(S_{n}) \ra C(S_{n-1}) \rtimes \zz$ be the embedding given by extension to zero outside $S_{n}$, composed with the canonical embedding $C(S_{n-1}) \subset C(S_{n-1} \rtimes \zz)$.

\begin{lem}
\label{lemin}
For every $n \in \nn$ there is a $*$-homomorphism $\phi_{n} : C(S_{n} )\rtimes_{\Phi_{n}} \zz \ra C(S_{n-1})\rtimes_{\Phi_{n-1}} \zz$ such that $\phi_{n} (f) =e_{n}(f)$ for every $f \in C(S_{n})$.
\end{lem}
\proof 

For $i=0,1$ denote by $u_{n-1+i}$ the unitary in $C(S_{n-1+i})\rtimes \zz$ implementing the action. Consider the Rokhlin towers decomposition $S_n = \bigsqcup_{j \in J_{n+1,n}}  \bigsqcup_{k=0}^j F_{n,j}^{(k)}$ relative to the inclusion $S_{n} \subset S_{n-1}$. The element $v_{n-1} := \sum_{i\in J_{n+1,n}} \chi_{\Phi_{n-1}^{j+1} (F_{n,j}^{(0)})} u_{n-1}^{i+1}$ is a partial isometry with support and range projections both equal to $\chi_{S_n}$. They give the covariance relation $v_{n-1} e_n (f) v_{n-1}^* = e_n (f \circ \Phi_n^{-1})$ for $f \in C(S_n)$. It follows that the map $\phi_n : C_c (\zz , C(S_n)) \ra C_c (\zz, C_c (S_{n-1}))$ given by $\phi_n (\sum_{k \in \zz} f_k u_n^k) = \sum_{k \geq 0} e_n(f) v_{n-1}^k + \sum_{k <0} f_k v_{n-1}^{k,*}$ is an algebraic $*$-homomorphism. Abusing notation, we write $v_{n-1}^* = v_{n-1}^{-1}$, in such a way that $\phi_n$ reads $\phi_n(\sum_{k \in \zz} f_k u_n^k)= \sum_{k \in \zz} e_n(f) v_{n-1}^k$. We want to check continuity of $\phi_n$. Let $x \in S_n$. This determines a sequence $\{m(j)\}_{j \in \zz} \subset \zz$ given by the condition $\Phi_{n-1}^m (x) \in S_n$ if and only if $m \in \{m(j)\}$. Consider the associated representation $\rho_x : C(S_{n-1} \rtimes \zz) \ra \mathbb{B}(l^2 (\zz))$. Let $\xi_x \in l^2 (\zz)$. Then
\begin{equation*}
\begin{split}
(\rho_x (\phi_n & (\sum_{k \in \zz} f_k u_n^k))\xi_x) (m) = \begin{cases}
\sum_{k \in \zz} (e_n (f) \circ \Phi_{n-1}^m)(x) \xi_x (m(j-k)) & \textit{ for } m=m(j),\\
0 & \textit{ otherwise }
\end{cases}\\
&= \begin{cases}
\sum_{k \in \zz} (e_n (f \circ \Phi_{n}^j))(x) \xi_x (m(j-k)) & \textit{ for } m=m(j),\\
0 & \textit{ otherwise }.
\end{cases}
\end{split}
\end{equation*}
It follows that $\phi_n$ is continuous. $\Box$

\begin{prop}
\label{propn3}
For every $n \in \mathbb{N}$ there is a commutative diagram
\begin{equation*}
\begin{array}{ccccccccccc} 0 &\ra &\mathbb{Z} &\overset{\eta_{n+1}}{\ra}& C(\Phi_n^{-1} (S_{n+1}), \mathbb{Z}) &\overset{\beta_{n+1}}{\ra}& K_0 (D_{n+1})& \overset{\tilde{\gamma}_{n+1,*}}{\ra}& K_0 (C(X) \rtimes \mathbb{R})& \ra &0\\
  & &\downarrow & &\downarrow & &\downarrow & &\downarrow & & \\
0 &\ra &\mathbb{Z} &\overset{\eta_{n+2}}{\ra}& C(\Phi_{n+1}^{-1} (S_{n+2}), \mathbb{Z}) &\overset{\beta_{n+2}}{\ra}& K_0 (D_{n+2})& \overset{\tilde{\gamma}_{n+2,*}}{\ra}& K_0 (C(X) \rtimes \mathbb{R})& \ra &0
\end{array}
\end{equation*}
with exact rows except in $C(\Phi_n^{-1} (S_{n+1}), \mathbb{Z})$ and $C(\Phi_{n+1}^{-1} (S_{n+2}), \mathbb{Z})$.
\end{prop}
\proof 
We have to define the vertical arrows (see Proposition \ref{propn2}) and the modified arrows $\tilde{\gamma}_{n+i,*}$, $i=1,2$. The map $\mathbb{Z} \ra \mathbb{Z}$ is the identity and the map $K_0 (D_{n+1}) \ra K_0 (D_{n+2})$ is $\iota_{n+1, *}$. 

For what concerns the second vertical arrow, let $f \in C(\Phi_n^{-1} (S_{n+1}), \mathbb{Z})$. We want to find $g \in C(\Phi_{n+1}^{-1} (S_{n+2}), \mathbb{Z})$ such that $\iota_{n+1,*} (g)-\iota_{n+1,*}(\hat{\Phi}_{n+1, \mathbb{K},*}(g)) = \iota_{n,n+2,*}(f) - \iota_{n,n+2,*} (\hat{\Phi}_{n,\mathbb{K},*} (f)) $, where we are considering $f$ respectively $g$ as functions defined on $S_n$ respectively $S_{n+1}$, after extension to zero outside their domains. Consider the tower decompositions $S_{n} = \bigsqcup_{j \in J_{n+1,n}}  T_{n+1}^{(j)}$ relative to $S_{n+1} \subset S_n$ and $S_{n+1} = \bigsqcup_{j \in S_{n+2,n+1}} T_{n+2}^{(j)}$ relative to $S_{n+2} \subset S_{n+1}$. Let $y \in F_{n+2,i}^{(0)} \subset S_{n+2}$ for some $i \in J_{n+2,n+1}$. Denote by $y_0=y, y_1,..., y_i$ the points in the forward orbit of $y$ which intersect $S_{n+1} \backslash S_{n+2}$ at times less than the time needed to $y$ in order to return to $S_{n+2}$ ($y_k = \Phi_{n+1}^k (y)$, $0 \leq k \leq i$). In the same way, every $y_k$ belongs to the $0$-th floor of a Rokhlin tower of height $l(k)$ associated to $\Phi_{n}$; for every $0\leq k \leq i$ and $0\leq j \leq l(k)$ we denote by $y_{k,j}$ the points in the forward orbit of $y_k=y_{k,0}$ which intersect $S_n \backslash S_{n+1}$ before returning to $S_{n+1}$ ($y_{k,j} = \Phi_n^j y_k$). Hence
\begin{equation*}
\iota_{n,n+2,*} f(y) =  \tr \left(\begin{array}{cccc} f(y_{0,l(0)})	&		&		&0	\\
								& f(y_{1,l(1)}) 	&		&	\\
								& 		& \ddots	&	\\
						0		&		&		& f(y_{i,l(i)})	\end{array}\right)
\end{equation*}
and
\begin{equation*}
\iota_{n,n+2,*} \hat{\Phi}_{n,\mathbb{K},*} f(y) =  \tr \left(\begin{array}{cccc} f(\Phi_n^{-1}y_{0})	&		&		&0	\\
								& f(y_{0,l(0)}) 	&		&	\\
								& 		& \ddots	&	\\
						0		&		&		& f(y_{i-1,l(i-1)})	\end{array}\right).
\end{equation*}
Thus
\begin{equation*}
\iota_{n,n+2,*} f(y) -\iota_{n,n+2,*} \hat{\Phi}_{n,\mathbb{K},*} f(y)= f(y_{i,l(i)})-f(\Phi_n^{-1} (y)).
\end{equation*}

Keep the above notation. Every point of $\Phi_{n+1}^{-1} (S_{n+2})$ is of the form $y_i$ for some $y \in S_{n+2}$ belonging to the $0$-th floor $F_{n+2,i}^{(0)}$ of a Rokhlin tower $T_{n+2}^{(i)}$. Define $g \in C(\Phi_{n+1}^{-1}(S_{n+2}),\mathbb{Z}) \subset C(S_{n+1},\mathbb{Z})$ as $g(y_i)=f(y_{i,l(i)})$. We have

\begin{equation*}
\begin{split}
\iota_{n+1,*} g(y) &=  \tr \left(\begin{array}{cccc} g(y_{0})	&		&		&0	\\
								& g(y_{1}) 	&		&	\\
								& 		& \ddots	&	\\
						0		&		&		& g(y_{i})	\end{array}\right)= \tr
						\left(\begin{array}{cccc} 0	&		&		&0	\\
								&\ddots 	&		&	\\
								& 		& 0	&	\\
						0		&		&		& g(y_i)	\end{array}\right)\\
						&=\tr
						 \left(\begin{array}{cccc} 0	&		&		&0	\\
								&\ddots 	&		&	\\
								& 		& 0	&	\\
						0		&		&		& f(y_{i,l(i)})\end{array}\right)
						\end{split}
\end{equation*}
and
\begin{equation*}
\begin{split}
\iota_{n+1,*} \hat{\Phi}_{n+1, \mathbb{K},*}g(y) &= \tr \left(\begin{array}{cccc} g(\Phi_{n+1}^{-1}(y_{0}))	&		&		&0	\\
								& g(\Phi_{n+1}^{-1}(y_{1})) 	&		&	\\
								& 		& \ddots	&	\\
						0		&		&		& g(\Phi_{n+1}^{-1}(y_{i}))	\end{array}\right)\\
						&= \tr \left(\begin{array}{cccc} g(\Phi_{n+1}^{-1}(y_0))	&		&		&0	\\
								&0 	&		&	\\
								& 		& \ddots	&	\\
						0		&		&		& 0	\end{array}\right)
						=\tr
						 \left(\begin{array}{cccc} 0	&		&		&0	\\
								&\ddots 	&		&	\\
								& 		& 0	&	\\
						0		&		&		& f(\Phi_n^{-1}(y_0))\end{array}\right)
						\end{split}
\end{equation*}
Hence
\begin{equation*}
\iota_{n+1,n+2,*} g(y) - \iota_{n+1,n+2,*} \hat{\Phi}_{n+1, \mathbb{K},*}g(y)=f(y_{i,l(i)})-f(\Phi_n^{-1} (y)).
\end{equation*}

We fix a $*$-isomorphism $C(X) \rtimes \rr \simeq (C(S_2) \rtimes_{\Phi_2} \zz)\otimes \mathbb{K}$ and then consider $K_0 (C(S_2) \rtimes_{\Phi_2} \zz)$ in place of $K_0 (C(X) \rtimes \rr)$. The rightmost vertical arrow is $\id_* : K_0 (C(S_2) \rtimes_{\Phi_2} \zz) \ra K_0 (C(S_2) \rtimes_{\Phi_2} \zz)$. Let $\tilde{\gamma}_{2,*} = \alpha_{2,*}$, where $\alpha_2$ is the canonical embedding $C(S_2) \subset C(S_2) \rtimes \zz$. Let $\tilde{\gamma}_{3,*} := e_{3,*}$ (cfr. Lemma \ref{lemin}) and for $n >3$, let $\tilde{\gamma}_{n,*} :=  \phi_{3,*} \circ ... \circ \phi_{n-1,*} \circ e_{n,*}$.

We have to check commutativity on the rightmost square. Let $n \geq 2$ and  $f \in C(S_{n+1}, \mathbb{Z})$; using the Rokhlin tower decompositon relative to $S_{n+2} \subset S_{n+1}$, we can write $f=\sum_{j \in J_{n+2,n+1}}  f_{j}$, where for every $j \in J_{n+2,n+1}$, $f_{j} = \sum_{k=0}^j f_{j}^{(k)}$ corresponds to the decomposition of the corresponding tower into $j+1$ floors. Then, taking classes in $K_0(C(S_n) \rtimes \zz)$, for every $j$ we obtain

\begin{equation*}
\begin{split}
\phi_{n+1,*}& e_{n+2,*}\iota_{n+1,*} f_{j}= \phi_{n+1,*} \left(\begin{array}{cccc} f_{j}^{(0)}&		&		&0	\\
								&  f_{j}^{(1)} \circ \Phi_{n+1}	&		&	\\
								& 		& \ddots	&	\\
						0		&		&		& f_{j}^{(j)}\circ \Phi_{n+1}^j	\end{array}\right)\\
					&=\phi_{n+1,*}\left(\begin{array}{cccc} f_{j}^{(0)}	&		&		&0	\\
								& f_{j}^{(1)} 	&		&	\\
								& 		& \ddots	&	\\
						0		&		&		& f_{j}^{(j)}	\end{array}\right) = \phi_{n+1,*} (f_{i,j}) = e_{n+1,*} (f_{i,j}).
						\end{split}
\end{equation*}
A similar procedure applies to the case $n=1$.

Note that the first raw for $n=1$ is exact in virtue of Proposition \ref{propn2}. Hence, in order to prove exactness for $n >1$, it is enough to show that the vertical maps $\delta_{n+1} : C(\Phi_n^{-1} (S_{n+1}), \mathbb{Z}) \ra C(\Phi_{n+1}^{-1} (S_{n+2}), \mathbb{Z})$ and $\iota_{n+1,*}: K_0 (D_{n+1}) \ra K_0 (D_{n+2})$ are surjective for every $n$. Surjectivity of the map $\iota_{n+1,*}$ was established in Proposition \ref{propn2}. We claim that $\delta_{n+1}$ is surjective as well. Indeed, consider the restriction to $\Phi_{n+1}^{-1} (S_{n+2})$ of the decomposition $S_{n+1}= \bigsqcup_{j \in J_{n+1,n}} F_{n+1,j}^{(0)}$ of $S_{n+1}$ in $0$-th floors of Rokhlin towers relative to the inclusion $S_{n+1} \subset S_n$: $\Phi_{n+1}^{-1} (S_{n+2}) = \bigsqcup_{j \in J_{n+1,n} }\hat{F}_{n+1,j}^{(0)}$, with $\hat{F}_{n+1,j}^{(0)} =F_{n+1,j}^{(0)} \cap \Phi_{n+1}^{-1} (S_{n+2})$. For $j \in J_{n+1,n}$ and $g \in C(\Phi_{n+1}^{-1} (S_{n+2}), \zz)$, we let $f|_{\Phi_n^j(\hat{F}_{n+1,j}^{(0)})} =g \circ \Phi_n^{-j} |_{\Phi_n^j(\hat{F}_{n+1,j}^{(0)})}$, where as usual we identify $g$ with its extension on $S_n$. Let also $f|_{\Phi_n^{-1} (S_{n+1}) \backslash \bigsqcup_{j \in J_{n+1,n}} \Phi_n^j(\hat{F}_{n,j}^{(0)})} = 0$. Then $\delta_{n+1} (f) =g$ and the result follows. $\Box$

For every $n \in \nn$ we can define a map $t_n : S_{n+1} \ra \Phi_n^{-1} (S_{n+1})$ by sending $y$ to $\Phi_n^j (y)$ for $y \in F_{n+1,j}^{(0)}$, with $j \in J_{n+1,n}$. This was implicitly used to define the maps $\delta_n$ in Proposition \ref{propn3}. Note that, for every $n \in \nn$, $\Phi_n \circ t_n = \Phi_{n+1}$.

\begin{lem}
\label{lemt}
For every $n, k \in \nn$ we have
\begin{equation*}
t_n \circ t_{n+1} \circ ... \circ t_{n+k} (\Phi_{n+k+1}^{-1} (S_{n+k+2})) = \Phi_n^{-1} (S_{n+k+2}).
\end{equation*}
\end{lem}
\proof 
For every $n \in \nn$, $\Phi_{n+1}$ factors as $\Phi_{n+1}=\Phi_n \circ t_n$. Hence, for every $n, k \in \nn$, $\Phi_{n+k+1} = \Phi_{n+k} \circ t_{n+k} =... = \Phi_n \circ t_n \circ ... \circ t_{n+k}$. The result follows. $\Box$

\begin{thm}
\label{propn4}
Let $\gamma_* = \lim_n \tilde{\gamma}_{n,*} : K_0 ((C(X) \rtimes \rr)_{\bar{x}}) \ra K_0 (C(X) \rtimes \rr)$. Then $\gamma_*$ is an isomorphism of ordered groups.
\end{thm}
\proof 
In virtue of Proposition \ref{propn3} we always have maps
\begin{equation*}
0 \ra \mathbb{Z} \overset{\eta}{\ra} \lim_n C(\Phi_n^{-1} (S_{n+1}), \mathbb{Z}) \overset{\beta}{\ra} K_0 ((C(X) \rtimes \rr)_{\bar{x}}) 	\overset{\gamma_*}{\ra} K_0 (C(X) \rtimes \rr) \ra 0
\end{equation*}
with exact rows except possibly at $\lim_n C(\Phi_n^{-1} (S_{n+1}), \mathbb{Z})$. We want to prove that $\eta$ is surjective. Denote by $\delta_{n+1} : C(\Phi_n^{-1} (S_{n+1}), \mathbb{Z}) \ra C(\Phi_{n+1}^{-1} (S_{n+2}), \mathbb{Z})$ the morphism defined in Proposition \ref{propn3}. Let $n \in \nn$ and $f \in C(\Phi_n^{-1}(S_{n+1}),\mathbb{Z})$ be given; we will see that there are $k \in \nn$ and $m \in \mathbb{Z}$ such that $\delta_{n+k} \circ \delta_{n+k-1} \circ ... \circ \delta_{n+1} (f)= \eta_{n+k+1} (m)$.\\
There is an open set $U \subset \Phi_n^{-1} (S_{n+1})$ containing the point $z= \bigcap_{m \in \nn} \Phi_n^{-1} (S_{n+m+1})$ on which $f$ is constant, $f(y) = m$ for some $m \in \mathbb{Z}$, for every $y \in U$. Let $k \in \nn$ be such that $\Phi_n^{-1} (S_{n+k+2}) \subset U$. It follows from Lemma \ref{lemt} that $\delta_{n+k+1} \circ \delta_{n+k} \circ ... \circ \delta_{n+1} (f) = \eta_{n+k+2} (m)$. Hence $\gamma_*$ is an isomorphism of groups. 

The inclusion $\alpha_n : C(S_n) \subset C(S_n) \rtimes_{\Phi_n} \mathbb{Z}$ induces a surjective group homomorphism at the level of $K_0$-groups, which sends the classes of projections in $C(S_n) \otimes \mathbb{K}$ surjectively on the classes of projections in $(C(S_n) \rtimes_{\Phi_n} \mathbb{Z}) \otimes \mathbb{K}$. By Proposition \ref{propn3} the same is true for $\gamma_*$. It follows that the map $\gamma_* : K_0 ((C(X)\rtimes \mathbb{R})_{\bar{x}}) \ra K_0 (C(X) \rtimes \mathbb{R})$ is an isomorphism of ordered groups, being surjective on the classes of projections and an isomorphism at the level of groups.   $\Box$

Let $n \in \nn$ and denote by $A_{\bar{x}}^{(n)}$ the orbit breaking subalgebra in the point $\bar{x} \in S_n$ associated to the minimal action $\Phi_n \in \Aut (C(S_n))$ (\cite{putnam}).

\begin{cor}
For every $n \in \nn$, $(C(X) \rtimes \rr)_{\bar{x}} \simeq A_{\bar{x}}^{(n)} \otimes \mathbb{K}$.
\end{cor}
\proof The result follows from Elliott classification theorem for AF-algebras (\cite{eaf}). $\Box$

\section{Appendix}
This appendix contains basic facts about flows on metric spaces that are used throughout the paper.

\begin{lem}
\label{lem.ind1}
Let $X$ be a locally compact metric space and $\phi : \rr \curvearrowright X$ a flow. Denote by $\mathcal{P}_c(X)$ the set of compact subsets of $X$, endowed with the Hausdorff distance $d_H$. If $t_1<t_2$ and $Y: [t_1,t_2] \ra \mathcal{P}_c(X)$, $t \mapsto Y_t$ is continuous, then
\begin{equation*}
\bigcup_{t \in [t_1,t_2]} (t, Y_t) \subset [t_1,t_2] \times X
\end{equation*}
is closed.
\end{lem}
\proof 
Let $(r,x) \in ([t_1,t_2] \times X )\backslash \bigcup_{t \in [t_1,t_2]} (t, Y_t)$. Let $\bar{d}$ be the metric on $[t_1,t_2] \times X$ given by $\bar{d} ((t_1,x_1), (t_2,x_2))= |t_1-t_2| + d(x_1,x_2)$. This metric generates the topology of $[t_1,t_2] \times X$. \\
Since $Y_r$ is closed, there is $\epsilon >0$ such that 
\begin{equation}
\label{eqdist}
\inf_{y_r \in Y_r} d(x, y_r)>\epsilon
\end{equation}
It follows from the hypothesis that there is $\epsilon /8 \geq \delta(\{ Y_t\})=\delta  >0$ such that
\begin{equation*}
d_H (Y_s, Y_t) < \epsilon /8 \qquad \mbox{ for every } |s-t|<  \delta.
\end{equation*}
In particular, using compactness of the $Y_t$'s, for every $y_s \in Y_s$ there is $\bar{y}_t \in Y_t$ such that
\begin{equation}
\label{eq.ind1}
d( y_s, \bar{y}_t) = \inf_{y_t \in Y_t} d( y_s, y_t ) \leq \sup_{y'_s \in Y_s} \inf_{y_t \in Y_t} d(y'_s, y_t) \leq d_H (Y_s, Y_t) < \epsilon /8
\end{equation}
for every $|s-t| < \delta$.\\
Let $\delta' ((r,x)) = \delta' >0$ be such that
\begin{equation*}
\bar{d} ((s,x'), (r,x)) < \epsilon /4 \qquad \mbox{ for every }|s-r| < \delta', \; x' \in B_{\delta'} (x).
\end{equation*}
Let $y_s \in Y_s$, $0<|s-r| < \min \{ \delta, \delta'\}$ and $x' \in B_{\delta'} (x)$. Let $\bar{y}_r$ be given as in (\ref{eq.ind1}). Then 
\begin{equation*}
\begin{split}
\epsilon &\overset{(\ref{eqdist})}{<} \bar{d} ((r,x), (r,\bar{y}_r)) \leq \bar{d} ((r,x), (s,x')) + \bar{d} ((s,x'), (s, y_s)) + \bar{d} ((s, y_s), (r,\bar{y}_r))\\
&\overset{(\ref{eq.ind1})}{<}\epsilon /4 + \bar{d} ((s,x'), (s, y_s)) +2( \epsilon /8)
\end{split}
\end{equation*}
Hence $\inf_{B_{\delta'} (x)}\{\bar{d}((s,x'), (s,y_s))\} > \epsilon /2$ for every $y_s \in Y_s$, with $|s-t| < \min \{ \delta, \delta'\}$ and so there is an open neighborhood of $(r,x)$ that is not contained in $\bigcup_{t \in [t_1,t_2]} (t, Y_t)$, proving that $\bigcup_{t \in [t_1,t_2]} (t, Y_t)$ is closed. $\Box$\\

\begin{lem}
\label{lem.ind2}
Let $X$ be a locally compact metric space and $\phi : \rr \curvearrowright X$ a flow. Let $Y \subset X$ compact, $\epsilon >0$ and $f \in C_0(X)$ be such that $|f||_Y \leq \epsilon$. There is $\delta >0$ such that $|f(x)| < 2\epsilon$ for every $x \in \phi_{[-\delta, +\delta]}(Y)$.
\end{lem}
\proof 
If the claim does not hold, there are sequences $\{y_n\} \subset Y$ and $t_n \ra 0$ such that $f(\phi_{t_n} (y_n)) \geq 2\epsilon$ for every $n \in \nn$. By compactness we can replace $\{y_n\}$ with a convergent subsequence $y_n \ra \bar{y} \in Y$. Then, by continuity of the flow, $\bar{y} = \lim \phi_{t_n}(y_n)$ and $f(\bar{y}) \geq 2\epsilon$, which is a contradiction. $\Box$\\

\begin{lem}
\label{lem.ind3}
Let $X$ be a locally compact metric space and $\phi : \rr \curvearrowright X$ a flow. For every $x \in X$, $t \in \rr$, $\epsilon >0$, there is $\delta >0$ such that $\phi_s (\mathcal{B}_{\delta} (x)) \subset \mathcal{B}_\epsilon (\phi_s (x))$ for every $ |s| \leq |t|$.
\end{lem}
\proof 
Suppose there is $\epsilon >0$ for which the claim does not hold. Then there are sequences $\delta_n \ra 0$, $y_n \in \mathcal{B}_{\delta_n} (x)$, $s_n \in [-|t|, + |t|]$ such that $d (\phi_{s_n} (y_n), \phi_{s_n} (x)) \geq \epsilon$ for every $n \in \nn$. Passing to a subsequence we can suppose $s_n \ra \bar{s} \in [-|t|, + |t|]$. By continuity of the flow $\epsilon \leq \lim d(\phi_{s_n} (y_n), \phi_{s_n} (x)) = d (\phi_{\bar{s}} (x), \phi_{\bar{s}} (x))=0$, a contradiction. $\Box$

\begin{lem}
\label{lemb}
Let $X$ be a locally compact metric space and $\phi : \rr \curvearrowright X$ be a flow. Suppose there are $x\in X$, $C \subset X$ closed, $t  \geq 0$ such that $C \cap \phi_{[0,t]} (x) = \emptyset$. Then there are non-empty open sets $U$ containing $x$ and $V$ containing $C$ such that $V \cap \phi_{[0,t]} (U ) = \emptyset$.
\end{lem}
\proof
Suppose this is not the case. Then for every choice of $U$ containing $x$ and $V$ containing $C$ there is $0 \leq s \leq t$ such that $\phi_s (U) \cap V \neq \emptyset$. In particular, we can take sequences of open sets $U^n \subset B_{1/n} (x)$ around $x$, $V_n (C)$ with $\cap_{n \in \nn} V_n (C) = C$ and $\{s_n \} \subset [0,t]$ such that $\phi_{s_n} (U^n) \cap V_n \neq \emptyset$ for every $n \in \nn$. Hence $U_n \cap \phi_{-s_n} (V_n) \neq \emptyset$ for every $n \in \nn$, $\cap_n (U_n \cap \phi_{-s_n} (V_n)) = \{ x\}$; thus $\phi_{s_n} (x) \in V_n$ for every $n \in\nn$. Passing to a subsequence we find $s \in [0,t]$ such that $\phi_s (x) \in C$, a contradiction. $\Box$\\

\begin{lem}
\label{lemaa}
Let $X$ be a locally compact metric space and $\phi : \rr \curvearrowright X$ a free flow. Suppose there is a flowbox $B$ with non-empty interior admitting a central slice $C$ which is a Cantor space and let $C' \subset C \cap B^\circ$ clopen. Then there is $r >0$ such that $\phi (C' \times (-r,+r))$ is open in $X$.
\end{lem}
\proof
Note that if $\phi (C' \times (-r,+r))$ is contained in $B^\circ$, then it is open in $X$, being the image of an open set under the map $\phi$. Hence it is enough to show that there is $r >0$ such that $\phi (C' \times (-r,+r)) \subset B^\circ$. Suppose this is not the case. Then there are sequences $\{y_n\} \subset C'$ and $r_n \ra 0$ such that $\phi (y_n , r_n) \notin B^\circ$ for every $n \in \nn$. By compactness of $C'$ we can suppose that $y_n \ra \bar{y} \in C'$. Hence $\phi (y_n , r_n) \ra \phi (\bar{y},0) \in C' \subset B^\circ$, contradicting the fact that $\phi (y_n , r_n) \notin B^\circ$ for every $n \in \nn$. $\Box$

 \section{Acknowledgments}
The author thanks the anonymous referee for the comments concerning this work which lead to an improved exposition. He also thanks Professor Wilhelm Winter for suggesting the topic of $C^*$-algebras associated to flows. This work is supported by the grant "Horizon 2020 - Quantum algebraic structures and models", CUP: E52I15000700002 and the MIUR - Excellence Departments - grant: "$C^*$-algebras associated to p-adic groups, bi-exactness and topological dynamics", CUP: E83C18000100006. It is part of the project "Interaction of Operator Algebras with Quantum Physics and Noncommutative Structure", supported by the grant "Beyond Borders", CUP: E84I19002200005. The author acknowledges the support of INdAM-GNAMPA.
 
\noindent \textbf{Conflict of interest} The author declares that he has no conflict of interest.

 

\end{document}